# Fuzziness and Funds Allocation in Portfolio Optimization

**Jack Allen[1], Sukanto Bhattacharya[2], Florentin Smarandache[3]**

[1] School of Accounting and Finance, Griffith University, Queensland, Australia

[2] School of Information Technology, Bond University, Queensland, Australia

[3] Department of Mathematics, University of New Mexico, Gallup, USA

**Abstract.**

Each individual investor is different, with different financial goals, different levels of risk tolerance and different personal preferences. From the point of view of investment management, these characteristics are often defined as objectives and constraints. Objectives can be the type of return being sought, while constraints include factors such as time horizon, how liquid the investor is, any personal tax situation and how risk is handled. It's really a balancing act between risk and return with each investor having unique requirements, as well as a unique financial outlook – essentially a constrained utility maximization objective. To analyze how well a customer fits into a particular investor class, one investment house has even designed a structured questionnaire with about two-dozen questions that each has to be answered with values from 1 to 5. The questions range from personal background (age, marital state, number of children, job type, education type, etc.) to what the customer expects from an investment (capital protection, tax shelter, liquid assets, etc.). A *fuzzy logic system* has been designed for the evaluation of the answers to the above questions. We have investigated the notion of fuzziness with respect to funds allocation.

**Introduction.**

In this paper we have designed our fuzzy system so that customers are classified to belong to any one of the following three categories: [1]

> **\*Conservative and security-oriented (risk shy)**
>
> **\*Growth-oriented and dynamic (risk neutral)**
>
> **\*Chance-oriented and progressive (risk happy)**

Besides being useful for clients, investor classification has benefits for the professional investment consultants as well. Most brokerage houses would value this information as it gives them a way of targeting clients with a range of financial products more effectively - including insurance, saving schemes, mutual funds, and so forth. Overall, many responsible brokerage houses realize that if they provide an effective service that is *tailored* to individual needs, in the long-term there is far more chance that they will retain their clients no matter whether the market is up or down.

Yet, though it may be true that investors can be categorized according to a limited number of types based on theories of personality already in the psychological profession's armory, it must be said that these classification systems based on the Behavioral Sciences are still very much in their infancy and they may still suffer from the problem of their meanings being similar to other related typographies, as well as of greatly oversimplifying the different investor behaviors. [2]

**(I.1) Exploring the implications of utility theory on investor classification.**

In our present work, we have used the familiar framework of neo-classical utility theory to try and devise a structured system for investor classification according to the utility preferences of individual investors (and also possible re-ordering of such preferences).

The theory of consumer behavior in modern microeconomics is entirely founded on observable utility preferences, rejecting hedonistic and introspective aspects of utility. According to modern utility theory, utility is a representation of a set of mutually consistent choices and not an explanation of a choice. The basic approach is to ask an individual to reveal his or her personal utility preference and not to elicit any numerical measure. [1] However, the projections of the consequences of the options that we face and the subsequent choices that we make are shaped by our memories of past experiences – that "mind's eye sees the future through the light filtered by the past". However, this memory often tends to be rather selective. [9] An investor who allocates a large portion of his or funds to the risky asset in period t-1 and makes a significant gain will perhaps be induced to put an even larger portion of the available funds in the risky asset in period t. So this investor may be said to have displayed a very weak risk-aversion attitude up to period t, his or her actions being mainly determined by past happenings one-period back. There are two interpretations of utility – normative and positive. *Normative utility* contends that *optimal* decisions do not always reflect the *best* decisions, as maximization of instant utility based on selective memory may not necessarily imply maximization of total utility. This is true in many cases, especially in the areas of health economics and social choice theory. However, since we will be applying utility theory to the very specific area of funds allocation between risky and risk-less investments (and investor classification based on such allocation), we will be concerned with *positive utility,* which considers the optimal decisions as they are, and not as what they should be. We are simply interested in using utility functions to classify an individual investor's attitude towards bearing risk at a given point of time. Given that the neo-classical utility preference approach is an objective one, we feel it is definitely more amenable to formal analysis for our purpose as compared to the philosophical conceptualizations of pure hedonism if we can accept decision utility preferences generated by selective memory.

If u is a given utility function and w is the wealth coefficient, then we have $E[u(w+k)] = u[w + E(k) – p]$, that is, $E[u(w+k)] = u(w - p)$, where k is the outcome of a risky venture given by a known probability distribution whose expected value $E(k)$ is zero. Since the outcome of the risky venture is as likely to be positive as negative, we would be willing to pay a small amount p, the *risk premium*, to avoid having to undertake the risky venture. Expanding the utilities in *Taylor series* to second order on the left-hand side and to first order on the right-hand side and subsequent algebraic simplification leads to the general formula $p = - (v/2) u''(w)/u'(w)$, where $v = E(k^2)$ is the variance of the possible outcomes. This shows that approximate risk premium is proportional to the variance – a notion that carries a similar implication in the mean-variance theorem of classical portfolio theory. The quantity $–u''(w)/u'(w)$ is termed the *absolute risk aversion*. [6] The nature of this absolute risk aversion depends on the form of a specific utility function. For instance, for a logarithmic utility function, the absolute risk aversion is dependent on the wealth coefficient w, such that it decreases with an increase in w. On the other hand, for an exponential utility function, the absolute risk aversion becomes a constant equal to the reciprocal of the risk premium.

**(I.2) The neo-classical utility maximization approach.**

In its simplest form, we may formally represent an individual investor's utility maximization goal as the following mathematical programming problem:

**Maximize U = f (x, y)**

**Subject to x + y = 1,**

**x ≥ 0 and y is unrestricted in sign**

Here x and y stand for the proportions of investable funds allocated by the investor to the market portfolio and a risk-free asset. The last constraint is to ensure that the investor can never borrow at the market rate to invest in the risk-free asset, as this is clearly unrealistic - the market rate being obviously higher than the risk-free rate. However, *an overtly aggressive investor can borrow at the risk-free rate* to invest in the market portfolio. In investment parlance this is known as ***leverage***. [5]

As in classical microeconomics, we may solve the above problem using the *Lagrangian multiplier* technique. The transformed Lagrangian function is as follows:

$$Z = f(x, y) + \lambda (1-x-y) \qquad (1)$$

By the first order (necessary) condition of maximization we derive the following system of linear algebraic equations:

$$Z_x = f_x - \lambda = 0 \quad \ldots(i)$$

$$Z_y = f_y - \lambda = 0 \quad \ldots(ii)$$

$$Z_\lambda = 1 - x - y = 0 \quad \ldots(iii) \qquad (2)$$

The investor's equilibrium is then obtained as the condition $f_x = f_y = \lambda^*$. $\lambda^*$ may be conventionally interpreted as the *marginal utility of money* (i.e. the investable funds at the disposal of the individual investor) when the investor's utility is maximized. [2]

The individual investor's indifference curve will be obtained as the locus of all combinations of x and y that will yield a constant level of utility. Mathematically stated, this simply boils down to the following total differential:

$$dU = f_x dx + f_y dy = 0 \qquad (3)$$

The immediate implication of (3) is that $dy/dx = -f_x/f_y$, i.e. assuming $(f_x, f_y) > 0$; this gives the negative slope of the individual investor's indifference curve and may be equivalently interpreted as the *marginal rate of substitution of allocable funds* between the market portfolio and the risk-free asset.

A second order (sufficient) condition for maximization of investor utility may be also derived on a similar line as that in economic theory of consumer behavior, using the sign of the bordered Hessian determinant, which is given as follows:

$$\overline{|H|} = 2\beta_x \beta_y f_{xy} - \beta_y^2 f_{xx} - \beta_x^2 f_{yy} \qquad (4)$$

In the above equation, $\beta_x$ and $\beta_y$ stand for the coefficients of x and y in the constraint equation. In this case we have $\beta_x = \beta_y = 1$. Equation (4) therefore reduces to:

$$\overline{|H|} = 2f_{xy} - f_{xx} - f_{yy} \qquad (5)$$

If $\overline{|H|} > 0$ then the stationary value of the utility function $U^*$ will be said to have attained its maximum.

To illustrate the application of classical utility theory in investor classification, let the utility function of a rational investor be represented by the following utility function:

$$U(x, y) = ax^2 - by^2; \text{ where}$$

x = proportion of funds invested in the market portfolio; and

y = proportion of funds invested in the risk-free asset.

Quite obviously, $x + y = 1$ since the efficient portfolio must consist of a combination of the market portfolio with the risk-free asset. The problem of funds allocation within the efficient portfolio then becomes that of maximizing the given utility function subject to the efficient portfolio constraint. As per J. Tobin's *Separation Theorem*; which states that investment is a two-phased process with the problem of portfolio selection which is considered independent of an individual investor's utility preferences (i.e. the first phase) to be treated *separately* from the problem of funds allocation within the selected portfolio which is dependent on the individual investor's utility function (i.e. the second phase). Using this concept we can mathematically categorize all individual investor attitudes towards bearing risk into any one of three distinct classes:

**Class A+: "Overtly Aggressive"(no risk aversion attitude)**

**Class A: "Aggressive" (weak risk aversion attitude)**
**Class B: "Neutral"(balanced risk aversion attitude)**
**Class C: "Conservative"(strong risk aversion attitude)**

The problem is then to find the general point of maximum investor utility and subsequently derive a mathematical basis to categorize the investors into one of the three classes depending upon the optimum values of x and y. The original problem can be stated as a *classical non-linear programming* with a single equality constraint as follows:

**Maximize $U(x, y) = ax^2 - by^2$**

**Subject to:**

**$x + y = 1$,**

**$x \geq 0$ and y is unrestricted in sign**

We set up the following transformed Lagrangian objective function:

**Maximize $Z = ax^2 - by^2 + \lambda(1 - x - y)$**

**Subject to:**

$x + y = 1$,

$x \geq 0$ and y is unrestricted in sign, (where $\lambda$ is the Lagrangian multiplier)

By the usual first-order (necessary) condition we therefore get the following system of linear algebraic equations:

$$Z_x = 2ax - \lambda = 0 \quad \text{... (i)},$$

$$Z_y = -2by - \lambda = 0 \quad \text{... (ii); and}$$

$$Z_\lambda = 1 - x - y = 0 \quad \text{... (iii)} \tag{6}$$

Solving the above system we get $x/y = -b/a$. But $x + y = 1$ as per the funds constraint. Therefore $(-b/a)y + y = 1$ i.e. $\mathbf{y^* = [1 + (-b/a)]^{-1} = [(a-b)/a]^{-1} = a/(a-b)}$. Now substituting for y in the constraint equation, we get $\mathbf{x^* = 1 - a/(a-b) = -b/(a-b)}$. Therefore the stationary value of the utility function is $\mathbf{U^* = a[-b/(a-b)]^2 - b[a/(a-b)]^2 = -ab/(a-b)}$.

Now, $f_{xx} = 2a$, $f_{xy} = f_{yx} = 0$ and $f_{yy} = -2b$. Therefore, by the second order (sufficient) condition, we have:

$$|\overline{H}| = 2f_{xy} - f_{xx} - f_{yy} = 0 - 2a - (-2b) = 2(b - a) \tag{7}$$

Therefore, the bordered Hessian determinant will be positive in this case if and only if we have $(a - b) < 0$. That is, given that $a < b$, our chosen utility function will be maximized at $U^* = ax^{*2} - by^{*2}$. However, the satisfaction of the non-negativity constraint on $x^*$ would require that $b > 0$ so that $-b < 0$; thus yielding $[-b/(a-b)] > 0$.

**Classification of investors:**

| Class | Basis of determination |
|-------|------------------------|
| A+    | $(y^* < x^*)$ and $(y^* \leq 0)$ |
| A     | $(y^* < x^*)$ and $(y^* > 0)$ |
| B     | $(y^* = x^*)$ |
| C     | $(y^* > x^*)$ |

**(I.3) Effect of a risk-free asset on investor utility.**

The possibility to lend or borrow money at a risk-free rate widens the range of investment options for an individual investor. The inclusion of the risk-free asset makes it possible for the investor to select a portfolio that *dominates* any other portfolio made up of only risky securities. This implies that an individual investor will be able to attain a higher *indifference curve* than would be possible in the absence of the risk-free asset. The risk-free asset makes it possible to *separate* the investor's decision-making process into two distinct phases – identifying the market portfolio and funds allocation. The market portfolio is the portfolio of risky assets that includes each and every available risky security. As all investors who hold any risky assets at all will choose to hold the market portfolio, this choice is independent of an individual investor's utility preferences.

Now, the expected return on a two-security portfolio involving a risk-free asset and the market portfolio is given by $E(R_p) = xE(R_m) + yR_f$; where $E(R_p)$ is the expected return on the optimal portfolio, $E(R_m)$ = expected return on the market portfolio; and $R_f$ is the return on the risk-free asset. Obviously, $x + y = 1$. Substituting for x and y with $x^*$ and $y^*$ from our illustrative case, we therefore get:

$$E(R_p)^* = [-b/(a-b)] E(R_m) + [a/(a-b)] R_f \qquad (8)$$

As may be verified intuitively, if b = 0 then of course we have E ($R_p$) = $R_f$, as in that case the optimal value of the utility function too is reduced to U* = -a0/(a-0) = 0.

The equation of the Capital Market Line in the original version of the CAPM may be recalled as **E ($R_p$) = $R_f$ + [E ($R_m$) - $R_f$]($S_p$/$S_m$)**; where E ($R_p$) is expected return on the efficient portfolio, E ($R_m$) is the expected return on the market portfolio, $R_f$ is the return on the risk-free asset, $S_m$ is the standard deviation of the market portfolio returns; and $S_p$ is the standard deviation of the efficient portfolio returns. Now, equating for E ($R_p$) with E ($R_p$)* we therefore get:

$$R_f + [E (R_m) - R_f](S_p/S_m) = [-b/(a-b)] \, E (R_m) + [a/(a-b)] \, R_f, \text{ i.e.}$$

$$S_p^* = S_m \, [R_f \{a/(a-b) – 1\} + \{-b/(a-b)\} \, E (R_m)] / [E (R_m) – R_f]$$

$$= S_m \, [E (R_m) – R_f][-b/(a-b)] / [E (R_m) – R_f]$$

$$= S_m \, [-b/(a – b)] \qquad (9)$$

This mathematically demonstrates that a rational investor having a quadratic utility function of the form U = $ax^2$ – $by^2$, at his or her point of maximum utility (i.e. affinity to return coupled with averseness to risk), assumes a given efficient portfolio risk (standard deviation of returns) equivalent to $S_p^*$ = $S_m$ [-b/(a – b)]; when the efficient portfolio consists of the market portfolio coupled with a risk-free asset.

The investor in this case, will be classified within a particular category A, B or C according to whether –b/(a-b) is greater than, equal in value or lesser than a/(a-b), given that a < b and b > 0.

**Case I: (b > a, b > 0 and a >0)**

Let b = 3 and a = 2. Thus, we have (b > a) and (-b < a). Then we have $x^* = -3/(2-3) = 3$ and $y^* = 2/(2-3) = -2$. Therefore $(x^* > y^*)$ and $(y^* < 0)$. So the investor can be classified as Class A+.

**Case II: (b > a, b > 0, a < 0 and b > |a|)**

Let b = 3 and a = - 2. Thus, we have (b > a) and (- b < a). Then, $x^* = -3/(-2-3) = 0.60$ and $y^* = -2/(-2-3) = 0.40$. Therefore $(x^* > y^*)$ and $(y^* > 0)$. So the investor can be re-classified as Class A!

**Case III: (b > a, b > 0, a < 0 and b = |a|)**

Let b = 3 and a = -3. Thus, we have (b > a) and (b = |a|). Then we have $x^* = -3/(-3-3) = 0.5$ and $y^* = -3/(-3-3) = 0.5$. Therefore we have $(x^* = y^*)$. So now the investor can be re-classified as Class B!

**Case IV: (b > a, b > 0, a < 0 and b < |a|)**

Let b = 3 and a = -5. Thus, we have (b>a) and (b<|a|). Then we have $x^* = -3/(-5-3) = 0.375$ and $y^* = -5/(-5-3) = 0.625$. Therefore we have $(x^* < y^*)$. So, now the investor can be re-classified as Class C!

So we may see that even for this relatively simple utility function, the final classification of the investor permanently into any one risk-class would be unrealistic as the range of values for the coefficients a and b could be switching dynamically from one range to another as the investor tries to adjust and re-adjust his or her risk-bearing attitude. This makes the neo-classical approach insufficient in itself to arrive at a classification. Here lies the justification to bring in a complimentary *fuzzy modeling* approach. Moreover, if we bring in time itself as an independent variable into the utility maximization framework, then one choice variable (weighted in favour of risk-avoidance) could be viewed as a *controlling factor* on the other choice variable (weighted in favour of risk-acceptance). Then the resulting problem could be gainfully explored in the light of *optimal control theory*.

**(II.1) Modeling fuzziness in the funds allocation behavior of an individual investor.**

The boundary between the preference sets of an individual investor, for funds allocation between a risk-free asset and the risky market portfolio, tends to be rather fuzzy as the investor continually evaluates and shifts his or her position; unless it is a passive ***buy-and-hold*** kind of portfolio.

Thus, if the universe of discourse is **U = {C, B, A and A+}** where C, B, A and A+ are our four risk classes "conservative", "neutral", "aggressive" and "overtly aggressive" respectively, then the fuzzy subset of U given by **P = {$x_1$/C, $x_2$/B, $x_3$/A, $x_4$/A+}** is the *true* preference set for our purposes; where we have **0 ≤ ($x_1$, $x_2$, $x_3$, $x_4$) ≤ 1**, all the symbols having their usual meanings. Although theoretically any of the P ($x_i$) values could be equal to unity, in reality it is far more likely that P ($x_i$) < 1 for i = 1, 2, 3, 4 i.e. the fuzzy subset P is most likely to be *subnormal*. Also, similarly, in most real-life cases it is expected that P ($x_i$) > 0 for i = 1, 2, 3, 4 i.e. all the elements of P will be included in its *support*: **supp (P) = {C, B, A, A+} = U**.

The critical point of analysis is definitely the individual investors preference ordering i.e. whether an investor is ***primarily conservative*** or ***primarily aggressive***. It is understandable that a primarily conservative investor could behave aggressively at times and vice versa but in general, their behavior will be in line with their classification. So the classification often depends on the height of the fuzzy subset P: **height (P) = $Max_x$P (x)**.

So one would think that the risk-neutral class becomes largely superfluous, as investors in general will tend to get classified as either primarily conservative or primarily aggressive. However, as already said, in reality, the element B will also generally have a *non-zero degree of membership* in the fuzzy subset and hence cannot be dropped.

The fuzziness surrounding investor classification stems from the fuzziness in the preference relations regarding the allocation of funds between the risk-free and the risky assets in the optimal portfolio. It may be mathematically described as follows:

Let **M** be the set of allocation options open to the investor. Then, the fuzzy preference relation is a fuzzy subset of the **M x M** space identifiable by the following membership function:

$$\mu_R (m_i, m_j) = 1; \quad m_i \text{ is definitely preferred to } m_j$$
$$c \in (0.5, 1); \quad m_i \text{ is somewhat preferred to } m_j$$
$$0.5; \quad \text{point of perfect neutrality}$$
$$d \in (1, 0.5); \quad m_j \text{ is somewhat preferred to } m_i; \text{ and}$$
$$0; \quad m_j \text{ is definitely preferred to } m_i \quad (10)$$

The fuzzy preference relation is assumed to meet the necessary conditions of reciprocity and transitivity. However, owing to substantial confusion regarding acceptable working definition of transitivity in a fuzzy set-up, it is often entirely neglected thereby leaving only the reciprocity property. This property may be succinctly represented as follows:

$$\mu_R (m_i, m_j) = 1 - \mu_R (m_j, m_i), \forall i \neq j \quad (11)$$

If we are to further assume a reasonable cardinality of the set **M**, then the preference relation $R_v$ of an individual investor v may also be written in a matrix form as follows: [12]

$$[r_{ij}^v] = [\mu_R (m_i, m_j)], \forall i, j, v \quad (12)$$

Classically, given the efficient frontier and the risk-free asset, there can be one and only one optimal portfolio corresponding to the point of tangency between the risk-free rate and the convex efficient frontier. Then fuzzy logic modeling framework does not in any way disturbs this bit of the classical framework. The fuzzy modeling, like the classical Lagrangian multiplier method, comes in only *after* the optimal portfolio has been identified and the problem facing the investor is that of allocating the available funds between the risky and the risk-free assets subject to a governing budget constraint. The investor is theoretically faced with an infinite number of possible combinations of the risk-free asset and the market portfolio but the ultimate allocation depends on the investor's utility function to which we now extend the fuzzy preference relation.

The available choices to the investor given his or her utility preferences determine the universe of discourse. The more uncertain are the investor's utility preferences, the wider is the range of available choices and the greater is the degree fuzziness involved in the preference relation, which would then extend to the investor classification. Also, wider the range of available choices to the investor the higher is the expected information content or *entropy* of the allocation decision.

**(II.2) Entropy as a measure of fuzziness.**

The term entropy arises in analogy with *thermodynamics* where the defining expression has the following mathematical form:

$$S = k \log_b \omega \qquad (13)$$

In thermodynamics, entropy is related to the *degree of disorder* or configuration probability $\omega$ of the canonical assembly. Its use involves an analysis of the microstates' distribution in the canonical assembly among the available energy levels for both isothermal reversible and isothermal irreversible (spontaneous) processes (with an attending modification). The physical scale factor k is the *Boltzmann constant*. [7]

However, the thermodynamic form has a different sign and the word *negentropy* is therefore sometimes used to denote expected information. Though Claude Shannon originally conceptualized the entropy measure of expected information, it was DeLuca and Termini who brought this concept in the realms of fuzzy mathematics when they sought to derive a universal mathematical measure of fuzziness.

Let us consider the fuzzy subset $F = \{r_1/X, r_2/Y\}$, $0 \leq (r_1, r_2) \leq 1$, where X is the event (y<x) and Y is the event (y≥x), x being the proportion of funds to be invested in the market portfolio and y being the proportion of funds to be invested in the risk-less security. Then the **DeLuca-Termini conditions** for measure of fuzziness may be stated as follows: [3]

- **FUZ (F) = 0** if F is a crisp set i.e. if the investor classified under a particular risk category *always* invests entire funds either in the risk-free asset (conservative attitude) or in the market portfolio (aggressive attitude)
- **FUZ (F) = Max FUZ (F)** when **F = (0.5/X, 0.5/Y)**
- **FUZ (F) ≥ FUZ (F*)** if F* is a *sharpened version* of F, i.e. if F* is a fuzzy subset satisfying $F^*(r_i) \geq F(r_i)$ given that $F(r_i) \geq 0.5$ and $F(r_i) \geq F^*(r_i)$ given that $0.5 \geq F(r_i)$

The second condition is directly derived from the concept of entropy. Shannon's measure of entropy for an n – events case is given as follows: [10]

$$H = -k \sum (p_i \log p_i), \text{ where we have } \sum p_i = 1 \qquad (14)$$

The Lagrangian form of the above function is as follows:

$$H_L = -k \sum(p_i \log p_i) + \lambda (1 - \sum p_i) \qquad (15)$$

Taking partial derivatives w.r.t. $p_i$ and setting equal to zero as per the necessary condition of maximization, we have the following stationary condition:

$$\frac{\partial H_L}{\partial p_i} = -k[\log p_i + 1] - \lambda = 0 \qquad (16)$$

It may be derived from (16) that at the point of maximum entropy, $\log p_i = -[(\lambda/k)+1]$, i.e. $\log p_i$ becomes a constant. This means that at the point of maximum entropy, $p_i$ becomes independent of the i and equalized to a constant value for i = 1, 2, ..., n. In an n-events case therefore, at the point of maximum entropy we necessarily have:

$$p_1 = p_2 = \ldots = p_i = \ldots = p_n = 1/n \qquad (17)$$

For n = 2 therefore, we obviously have the necessary condition for entropy maximization as $p_1 = p_2 = ½ = 0.5$. In terms of the Fuzzy preference relation, this boils down to exactly the second DeLuca-Termini condition. Keeping this close relation with mathematical information theory in mind, DeLuca and Termini even went on to incorporate Shannon's entropy measure as their chosen measure of fuzziness. For our portfolio funds allocation model, this measure could simply be stated as follows:

$$\text{FUZ}(F) = -k[\{F(r_1) \log F(r_i) + (1-F(r_1)) \log (1-F(r_1))\} + \{F(r_2) \log F(r_2) + (1-F(r_2)) \log (1-F(r_2))\}] \qquad (18)$$

**(II.3) Metric measures of fuzziness.**

Perhaps the best method of measuring fuzziness will be through measurement of the distance between F and $F^c$, as fuzziness is mathematically equivalent to the lack of distinction between a set and its complement. In terms of our portfolio funds allocation model, this is equivalent to the ambivalence in the mind of the individual investor regarding whether to put a larger or smaller proportion of available funds in the risk-less asset. The higher this ambivalence, the closer F is to $F^c$ and greater is the fuzziness.

This measure may be constructed for our case by considering the fuzzy subset F as a *vector with 2 components*. That is, F ($r_i$) is the i$^{th}$ component of a vector representing the fuzzy subset F and (1 – F ($r_i$)) is the i$^{th}$ component of a vector representing the complementary fuzzy subset $F^c$. Thus letting D be a metric in *2* space; we have the distance between F and $F^c$ as follows: [11]

$$D_\rho (F, F^c) = [\Sigma |F (r_i) - F^c (r_i)|^\rho]^{1/\rho}, \text{ where } \rho = 1, 2, 3, \ldots \quad (19)$$

For Euclidean Space with $\rho=2$, this metric becomes very similar to the statistical variance measure RMSD (root mean square deviation). Moreover, as $F^c (r_i) = 1 - F (r_i)$, the above formula may be written in a simplified manner as follows:

$$D_\rho (F, F^c) = [\Sigma |2F (r_i) - 1|^\rho]^{1/\rho}, \text{ where } \rho = 1, 2, 3, \ldots \quad (20)$$

For $\rho=1$, this becomes the *Hamming metric* having the following form:

$$D_1 (F, F^c) = \Sigma |2F(r_i) - 1| \quad (21)$$

If the investor always puts a greater proportion of funds in either the risk-free asset or the market portfolio, then F is reduced to a crisp set and $|2F (r_i) - 1| = 1$.

Based on the above metrics, a universal measure of fuzziness may now be defined as follows for our portfolio funds allocation model. This is done as follows:

For a crisp set F, $F^c$ is truly complementary, meaning that the metric distance becomes:

$$D_\rho *(F, F^c) = 2^{1/\rho}, \text{ where } \rho = 1, 2, 3, \ldots \quad (22)$$

An effective measure of fuzziness could therefore be as follows:

$$FUZ_\rho (F) = [2^{1/\rho} - D_\rho (F, F^c)] / 2^{1/\rho} = 1 - D_\rho (F, F^c) / 2^{1/\rho} \quad (23)$$

For the Euclidean metric we would then have:

$$FUZ_2(F) = 1 - \frac{[\Sigma(2F(r_i) - 1)^2]^{1/2}}{\sqrt{2}}$$

$$= 1 - (\sqrt{2})(RMSD), \text{ where } RMSD = ([\Sigma(2F(r_i) - 1)^2]^{1/2})/2 \qquad (24)$$

For the Hamming metric, the formula will simply be as follows:

$$FUZ_1(F) = 1 - \frac{\Sigma|2F(r_i) - 1|}{2} \qquad (25)$$

Having worked on the applicable measure for the degree of fuzziness of our governing preference relation, we devote the next section of our present paper to the incorporation of neurofuzzy control systems to produce and compute risk classification for the homological utility of an investor under different scenarios. In a subsequent section, we also take a passing look at the possible application of optimal control theory to model the dynamics of funds allocation behavior of an individual investor.

**(IV) Exploring time-dependent funds allocation behavior of individual investor in the light of optimal control theory.**

If the inter-temporal utility of an individual viewed from time t is recursively defined as $U_t = W[c_t, \mu(U_{t+1}| I_t)]$, then the aggregator function W makes current inter-temporal utility a function of current consumption $c_t$ and of a certainty equivalent of next period's random utility $I_t$ that is computed using information up to t. Then, the individual could choose a control variable $x_t$ in period t to maximize $U_t$.[4] In the context of the mean-variance model, a suitable candidate for the control variable could well be the proportion of funds set aside for investment in the risk-free asset. So, the objective function would incorporate the investor's total temporal utility in a given time range [0, T]. Given that we include time as a continuous variable in the model, we may effectively formulate the

problem applying classical optimal control theory. The plausible methodology for formulating this model is what we shall explore in this section.

The basic optimal control problem can be stated as follows: [8]

Find the control vector **u** = ($u_1$, $u_2$, …, $u_m$) which optimizes the functional, called the *performance index,* J = $\int f_0$ (**x**, **u**, t) dt over the range (0, T), where **x** = ($x_1$, $x_2$, …, $x_n$) is called the state vector, t is the time parameter, T is the terminal time and $f_0$ is a specified function of **x**, **u** and t. The state variables $x_i$ and the control variables $u_i$ are related as $dx_i/dt$ = $f_i$ ($x_1$, $x_2$, …, $x_n$; $u_1$, $u_2$, …, $u_m$; t), i = 1, 2, …, n.

In many control problems, the system is linearly expressible as **x** (.) = $[A]_{nxn}$ **x** + $[B]_{nxm}$ **u**, where all the symbols have their usual connotations. As an illustrative example, we may again consider the quadratic function that we used earlier $f_0$ (x, y) = $ax^2 - by^2$. Then the problem is to find the control vector that makes the performance index J = $\int(ax^2 - by^2)$ dt stationary with x = 1 − y in the range (0, T).

The ***Hamiltonian*** may be expressed as **H = $f_0$ + λy = ($ax^2 - by^2$) + λy**. The standard solution technique yields $-H_x$ = λ(.) … (i) and $H_u$ = 0 … (ii) whereby we have the following system of equations: -2ax = λ(.) … (iii) and -2y + λ = 0 … (iv). Differentiation of (iv) leads to –2y(.) + λ(.) = 0 … (v). Solving (iii) and (v) simultaneously, we get 2ax = -2y(.) = -λ(.) i.e. y(.) = -ax … (vi). Transforming (iii) in terms of x and solving the resulting ordinary differential equation would yield the state trajectory x (t) and the optimal control u (t) for the specified quadratic utility function, which can be easily done by most standard mathematical computing software packages.

So, given a particular form of a utility function, we can trace the dynamic time-path of an individual investor's fund allocation behavior (and hence; his or her classification) within the ambit of the mean-variance model by obtaining the state trajectory of x – the proportion of funds invested in the market portfolio and the corresponding control variable y – the proportion of funds invested in the risk-free asset using the standard techniques of optimal control theory.

*References:*

**Book, Journals, and Working Papers:**

[1] Arkes, Hal R., Connolly, Terry, and Hammond, Kenneth R., "Judgment And Decision Making – An Interdisciplinary Reader" Cambridge Series on Judgment and Decision Making, Cambridge University Press, 2$^{nd}$ Ed., 2000, 233-39

[2] Chiang, Alpha C., "Fundamental Methods of Mathematical Economics" McGraw-Hill International Editions, Economics Series, 3$^{rd}$ Ed. 1984, 401-408

[3] DeLuca A.; and Termini, S., "A definition of a non-probabilistic entropy in the setting of fuzzy sets", *Information and Control 20*, 1972, 301-12

[4] Haliassos, Michael and Hassapis, Christis, "Non-Expected Utility, Savings and Portfolios" *The Economic Journal,* Royal Economic Society, January 2001, 69-73

[5] Kolb, Robert W., "Investments" Kolb Publishing Co., U.S.A., 4$^{th}$ Ed. 1995, 253-59

[6] Korsan, Robert J., "Nothing Ventured, Nothing Gained: Modeling Venture Capital Decisions" *Decisions, Uncertainty and All That*, The Mathematica Journal, Miller Freeman Publications 1994, 74-80

[7] Oxford Science Publications, *Dictionary of Computing* OUP, N.Y. USA, 1984

[8] S.S. Rao, *Optimization theory and applications*, New Age International (P) Ltd., New Delhi, 2$^{nd}$ Ed., 1995, 676-80

[9] Sarin, Rakesh K. & Peter Wakker, "Benthamite Utility for Decision Making" Submitted Report, Medical Decision Making Unit, Leiden University Medical Center, The Netherlands, 1997, 3-20

[10] Swarup, Kanti, Gupta, P. K.; and Mohan, M., *Tracts in Operations Research* Sultan Chand & Sons, New Delhi, 8$^{th}$ Ed., 1997, 659-92

[11] Yager, Ronald R.; and Filev, Dimitar P., *Essentials of Fuzzy Modeling and Control* John Wiley & Sons, Inc. USA 1994, 7-22

[12] Zadrozny, Slawmir "An Approach to the Consensus Reaching Support in Fuzzy Environment", *Consensus Under Fuzziness* edited by Kacprzyk, Januz, Hannu, Nurmi and Fedrizzi, Mario, International Series in Intelligent Systems, USA, 1997, 87-90

**Website References:**

*1* http://www.fuzzytech.com/e/e_ft4bf6.html

*2* http://www.geocities.com/wallstreet/bureau/3486/11.htm